\documentclass{amsart}

\usepackage{rlepsf, boxedminipage, amssymb}

\newtheorem{thm}{Theorem}[section]
\newtheorem{lem}[thm]{Lemma}

\newtheorem{rem}[thm]{Remark}

\newtheoremstyle{definition}{7pt plus6.3pt minus6.3pt}{7pt plus3pt minus3pt}%
{\rm}{}{\bf}{}{0.75em}{\thmname{#1}\thmnumber{ #2}\thmnote{\sl\stdspace#3}}
\theoremstyle{definition}\newtheorem{example}[thm]{Example}
\newtheorem{exercise}[thm]{\small Exercise}

\newcommand{\bbr}{\begin{rem}\em}
\newcommand{\eer}{\end{rem}}

\newcommand{\krn}{\operatorname{ker}}
\newcommand{\la}{\langle}
\newcommand{\ra}{\rangle}

\newcommand{\bex}{\begin{example}}
\newcommand{\eex}{\end{example}}

\newcommand{\bhw}{\begin{exercise}\small}
\newcommand{\ehw}{\end{exercise}}

\newcommand{\pa}{\partial}
\newcommand{\be}{\begin{enumerate}}
\newcommand{\ee}{\end{enumerate}}

\def\CC{\hbox{$\mathcal{C}$ } }
\def\M{\hbox{$\mathcal{M}$ } }
\def\A{\hbox{$\mathcal{A}$ } }
\def\C{\hbox{$\mathbb C$} }
\def\Z{\hbox{$\mathbb Z$} }

\def\R{\hbox{$\mathbb R$} }

\def\co{\colon\thinspace}

\def\dfn#1{{\em #1}}

\begin{document}

\title{Invariants of Knots, Embeddings and Immersions via Contact Geometry}

\author{Tobias Ekholm}
\address{Department of Mathematics,
University of Southern California,
Los Angeles, CA 90089-1113}

\author{John B. Etnyre}
\address{University of Pennsylvania, Philadelphia, PA 19104}
\email{etnyre@math.upenn.edu}
\urladdr{http://www.math.upenn.edu/\char126 etnyre}


\begin{abstract}
This paper is an overview of the idea of using contact geometry to construct invariants of immersions
and embeddings. In particular, it discusses how to associate a contact manifold to any manifold and
a Legendrian submanifold to an embedding or immersion. We then discuss recent work that creates
invariants of immersions and embeddings using the Legendrian contact homology of the associated
Legendrian submanifold.
\end{abstract}

\maketitle

In recent years classical constructions in contact and symplectic geometry has been applied in the
study of topological problems. This idea was used by Arnold \cite{Aronld94} where he studies the
structure of the space of immersed plane curves. More specifically, Arnold uses a lift of an
immersed curve in the plane to the unit cotangent bundle of $\R^2$, which is diffeomorphic to
$\R^2\times S^1$ and which carries a natural contact structure. The lift is a Legendrian knot
provided the curve is self transverse and Arnold demonstrated that the contact invariants of this
Legendrian knot are useful for understanding the classification of plane curves up to regular
homotopies restricted in certain ways. Since then others have picked up on Arnold's perspective to
study codimension one immersions in other dimensions \cite{Goryunov}.

This had been an exciting area of research for some time when Ooguri and Vafa \cite{OV} suggested a
similar construction had relevance to modern physics (in particular ``large $N$-dualities'') a new
wave of interest in this construction began. In particular, the idea of using holomorphic curves (or
``branes'') on a Lagrangian submanifold to study knots in $\R^3$ seemed very promising. Moreover,
the development of Symplectic Field Theory by Eliashberg, Givental and Hofer \cite{egh} seemed like
the ideal theory to count these holomorphic curves. Recently \cite{ees3} has provided the analytic
underpinnings of the small part of this theory, called Legendrian contact homology in $1$-jet
spaces, that is necessary to define invariants of knots. The construction of Ooguri and Vafa has
inspired Ng \cite{Ng1, Ng2, Ng3} to define an invariant of knots using a braid description of the
knot. This invariant should compute the contact homology of the Legendrian submanifold associated to
a knot in $\R^3,$ though it is ongoing research to prove this. Non the less, Ng's invariant has been
shown to be amazingly powerful \cite{Ng3}. For example, it seems to contain the Alexander and $A$
polynomials, can distinguish the unknot from other knots and at the moment it seems possible that it
is a complete invariant of knots.

The goal of this paper is to give a somewhat topologist friendly introduction to this theory and
describe the connections between the invariant described by Ng and Legendrian contact homology. The
story involves a beautiful interplay between topology, geometry, analysis, combinatorics and
algebra. We will only begin the story and hopefully provide the background necessary for the reader
to pick up \cite{Ng3} or \cite{ees3} and learn more about it. In particular the two papers
concentrate on different parts of the story (the combinatorics and algebra in the fist paper and the
analysis and geometry in the second) and we hope the reader of this paper will see them both as part
of a unified theory that has promise far beyond what is currently in the literature.

Acknowledgments: We would like to thank Hans Boden for encouraging us to write this overview of the
material contained in a talk the second author gave at ``Geometry and Topology of Manifolds''
conference held at McMaster University in May of 2004. TE is a research fellow of the Swedish Royal 
Academy of Sciences sponsored by the Knut and Alice Wallenberg foundation.
JE was supported in part by
NSF CAREER Grant (DMS--0239600) and FRG-0244663.

\section{A Geometric Construction}

To a manifold $M$ we associate contact manifold $(W_M, \xi_M)$ and to an embedding or immersion
$N\subset M$ we associate a Legendrian submanifold $L_N$ in $(W_M,\xi_M).$ In the first subsection
we briefly recall the definition of a contact structure and a Legendrian submanifold. In
Subsection~\ref{ctmfd} we describe the contact manifold $(W_M, \xi_M)$ associated to $M$ and in the
last two sections we describe the Legendrian $L_N$ associated to an embedded, respectively immersed,
submanifold $N\subset M.$ To the reader not familiar with contact geometry we recommend consulting
the papers \cite{EtnyreCIntro, GeigesIntro}. It also might be useful to look at the following
section for specific examples of contact manifolds and Legendrian submanifolds.


\subsection{Basic contact notions}
An oriented {\em contact structure} on an orientable $(2m+1)$-manifold $M$ is a completely
non-integrable field of tangent hyperplanes $\xi\subset TM$. That is, a field of hyperplanes given
as $\xi=\krn{\alpha}$, where the non-vanishing $1$-form $\alpha$ (the {\em contact form}) is such
that the $(2m+1)$-form $\alpha\wedge(d\alpha)^m$ is a volume form on $M$. Note that if $\alpha$ is a
contact form then $d\alpha|\xi$ is a symplectic form. A diffeomorphism of contact manifolds
$(M,\xi)\to (M',\xi')$ is called a {\em contactomorphism} if it maps $\xi$ to $\xi'$. An immersion
of an $m$-manifold $f\colon L\to M$ is called {\em Legendrian} if $df_p(T_p L)\subset \xi_{f(p)}$
for all $p\in L$.

\subsection{A natural contact manifold}\label{ctmfd}
We begin with a manifold $M$ of dimension $m.$ The goal of this subsection is to put a contact
structure on the oriented projectivized co-tangent bundle $W_M.$ That is $W_M$ is the space of
oriented lines in $T^*M$:
\[
W_M=\{v\in T^*_x M: v\not=0\}/\sim,
\]
where $v\sim v'$ if $v=cv$ for a positive constant $c.$ Clearly $W_M$ is a $S^{m-1}$ bundle over
$M.$ It is sometimes helpful to think of $W_M$ using a metric. Fix a metric $g$ on $M$ then set
\[
W_g=\{v\in T^*_xM: |v|_g=1\}
\]
to be the unit co-tangent bundle of $M.$ Now there is an obvious map
\[\phi_g:W_M\to W_g\]
sending the equivalence class of $v$ to $\frac{v}{|v|_g}.$ It is clear that $\phi_g$ is a
diffeomorphism. We will usually call $W_M$ the unit co-tangent bundle, but it is sometime useful to
realize that $W_M$ can be defined without choosing a metric.

Recall there is a canonical 1-form $\lambda,$ usually called the Liouville 1-form, on $T^*M.$ The
1-form $\lambda$ is characterized by the property that for any 1-form $\alpha:M\to T^*M$ on $M$ the
pull-back of $\lambda$ to $M$ by $\alpha$ is $\alpha$:
\begin{equation}\label{charlam}
\alpha^* \lambda = \alpha
\end{equation}
We now express $\lambda$ in local coordinates. Let $q_1,\ldots, q_m$ be local coordinates in the
open set $U\subset M.$ Then coordinates on $T^*U\subset T^*M$ are $q_1,\ldots, q_m,p_1, \ldots, p_m$
where any point $\beta$ in $T^*_{(q_1,\ldots, q_m)}U$ can be written
\[
\beta =\sum_{i=1}^m p_i\, dq_i.
\]
Set $dq_i'=\pi^* dq_i,$ where $\pi\colon T^*U \to U$ is the projection map, then we can write
\[
\lambda = \sum_{i=1}^m p_i\, dq_i'.
\]
It is easy to verify that the right hand side satisfies \eqref{charlam}. It is customary to abuse
notation and write $dq_i$ for $dq_i'$ and let the context define the meaning of $dq_i.$ Though this
can sometimes be confusing we adopt this standard abuse of terminology here.

\begin{lem}
Thinking of $W_M$ as the unit cotangent bundle (having fixed some metric) then
\[
\alpha = \lambda|_{W_M}
\]
is a contact form on $W_M.$
\end{lem}
\begin{proof}
It is clear that
\[
d\lambda = \sum_{i=1}^m dp_i\wedge dp_i
\]
is a symplectic form on $T^*M.$ In particular $(d\lambda)^m\not=0.$ Note the vector $v=\sum_{i=1}^m
p_i \frac{\partial}{\partial p_i}$ is a vector field on $T^*M$ that is transverse to $W_M.$
Moreover, contracting $v$ into $d\lambda$ yields $\lambda.$ So
\[
\lambda \wedge (d\lambda)^{m-1}=\frac 1m \iota_v(d\lambda)^m
\]
is a $2m-1$ from that is not zero on any hyperplane in $T^*M$ that is transverse to $v.$ (Here
$\iota$ denotes contraction.) Thus it is not zero on the tangent planes to $W_M.$
\end{proof}
Thus $\xi=\ker \lambda|_{W_M}$ is a contact structure on $W_M.$

Just as we can define the manifold $W_M$ in a metric independent way, we can define the contact
structure $\xi$ without reference to a choice of metric on $M$. To this end note that we can
identify $W_M$ with the bundle of oriented tangent hyperplanes to $M$. This space is naturally
identified with the quotient space of the space of non-zero cotangent vectors along $M$ under the
identifictaion
$$
w\simeq a w,\quad a\in(0,\infty).
$$
The natural identification sends a covector to its kernel and its inverse sends a hyperplane to a
covector whose kernel is that hyperplane. We may coordinatize this quotient by using fiberwise
projective coordinates: if $(q,p)$ are the local coordinates in the proof above we may use
coordinates
$$
(q,[p_1,\dots, p_{k-1},\pm 1,p_{k+1},\dots,p_m]),\quad k=1,\dots,m,
$$
with the obvious transition functions. Using these coordinate patches we can think of the tangent
planes to the quotient space as affine cotangent hyperplanes and thus we may again restrict
$\lambda_M$. A straightforward calculation shows that in coordinates $(q,[p_1,\dots,p_{k-1},\pm
1,p_{k+1},\dots,p_m])$ we have
$$
\lambda\wedge(d\lambda)^{m-1}=\pm dq_k\wedge dp_1\wedge
dq_1\wedge\dots dp_{k-1}\wedge dq_{k-1}\wedge dp_{k+1}\wedge
dq_{k+1}\dots\wedge dq_m\wedge dp_m.
$$

\begin{center}
\begin{boxedminipage}[l]{3in}
\begin{center}
       To a {\bf smooth manifold} $M$ \break we associate the \break {\bf contact manifold}
       $(W_M,\xi).$
\end{center}
\end{boxedminipage}
\end{center}

\subsection{Submanifolds}
Now suppose we have a submanifold $N$ in $M.$ Let $L_N$ be the unit co-normal bundle
\[
L_N= \{u\in W_M | u(v)=0 \text{ for all } v \in TN\}.
\]
Clearly $L_N$ is a $(m-n-1)$-sphere bundle over $N$ and thus $L_N$ is a $(m-1)$-manifold in
$W_M^{2m-1}$ Suppose $u\in L_N$ and $v\in T_uL_N\subset T_uW_M.$ If $\pi$ denotes projection from
the tangent bundle to the base manifold then clearly $\pi_*(v)\in T_{\pi(u)}L.$ Hence
\[
\lambda_u(v)=u(\pi_*(v))=0
\]
by the definition of $L_N.$ So the contact form on $W_M$ vanishes on the tangent space to $L_N,$ or
in other words $T_uL_N\subset \xi_u$ for all $u\in L_N.$ This and the fact that $L_N$ is always
$m-1$ dimensional says that $L_N$ is a Legendrian submanifold of $(W_M,\xi).$
\begin{center}
\begin{boxedminipage}[l]{3in}
\begin{center}
  To a {\bf smooth submanifold} $N$ of $M$ \break we associate the \break {\bf Legendrian
  submanifold} $L_N$ of $(W_M,\xi).$
\end{center}
\end{boxedminipage}
\end{center}
Moreover, it is important to notice that as we isotop $N$ the Legendrian $L_N$ goes through a
Legendrian isotopy. Thus the Legendrian isotopy type of $L_N$ is an invariant of $N$ up to isotopy.

\subsection{Immersions}
Noting that the construction in the previous subsection is purely local we conclude that it works
more generally for immersed submanifold $N$ in $M$. In general, the result $L_N$ of the construction
is an immersed Legendrian submanifold of $W_M.$ Since Legendrian immersions satisfy an h-principle
there is little hope to derive subtle geometric invariants of $N$ from the Legendrian regular
homotopy class of $L_N$. However, $L_N$ is embedded provided that at any double point $p$ of $N$ in
$M$ there is no covector which annihilates both of the tangent space of $N$ at $p$. In particular,
$L_N$ is embedded for self transverse immersions.

If case $N$ a co-dimension 1 immersion then $L_N$ is the orientation double cover of $N$. In
particular, if $N$ is also oriented one may use ``half'' of the conormal lift: $L_N$ has two
components and the orientation picks out one preferred component of $L_N.$ Note that each component
of $L_N$ will be embedded as long as $N$ does not have double points where the tangent planes to $N$
agree as oriented hyperplanes. We will abuse notation and use $L_N$ to denote the preferred
component of $L_N$ whenever $N$ is an immersed oriented co-dimension 1 submanifold of $M.$ Thus
\begin{center}
\begin{boxedminipage}[l]{4in}
\begin{center}
       To a {\bf generically immersed submanifold} $N$ of $M$ \break we associate
       the \break {\bf Legendrian submanifold} $L_N$ of $(W_M,\xi).$
\end{center}
\end{boxedminipage}
\end{center}
The Legendrian isotopy class of $L_N$ is an invariant of $N$ up to regular homotopies avoiding
double points with tangencies as described above. That is whenever the regular homotopy goes through
a self tangency the orientations on the tangent planes must be opposite. Self tangencies where the
orientations on the tangent planes agree are called {\em dangerous self-tangencies}.

\section{Examples}
In this section we consider important examples of the construction outlined in the previous section.

\subsection{Plane Curves}
\begin{figure}[ht]
  \relabelbox \small {\epsfxsize=2in\centerline{\epsfbox{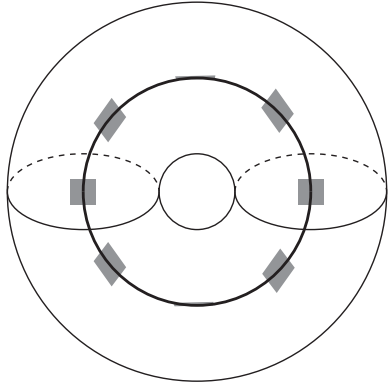}}}
  \endrelabelbox
        \caption{The unit cotangent bundle of $\R^2$. The thicker circle in the middle is the unit
          circle above the origin. A few contact planes along this circle are indicated. They make
          a complete turn as the circle is traversed.}
        \label{fig:UR2}
\end{figure}

The simplest example to consider is immersion of (oriented) $S^1$ into $\R^2.$ The unit cotangent
bundle is $W=\R^2\times S^1$ and the contact structure is $\alpha= (p_1\, dq_1 + p_2\,
dq_2)|_{W_M}.$ Using coordinates $(x_1, x_2, \theta)$ on $W,$ where $x_i=q_i$ and $\tan\theta=
\frac{p_1}{p_2},$ we can rewrite
\begin{align*}
\alpha&= p_1\, dq_1 + p_2\, dq_2 = p_2(\frac{p_1}{p_2}\, dq_1 + dq_2)\\ &= p_2 (\tan \theta \, dx_1
+ dx_2)= \frac{p_2}{\cos \theta}(\sin\theta\, dx_1 + \cos\theta \, dx_2)\\ &= \sin\theta\, dx_1 +
\cos\theta \, dx_2.
\end{align*}
The last inequality follows from $\frac{p_2}{\cos \theta}=\frac{p_2}{\frac{p_2}{p_1^2+p_2^2}}=
p_1^2+p_2^2=1.$ Thus the contact structure on $W$ is
\[
\xi=\ker(\sin\theta\, dx_1 + \cos\theta \, dx_2).
\]
This is easy to picture. See Figure~\ref{fig:UR2}.
All of the contact planes contain $\frac{\partial}{\partial \theta},$ so they are spanned by this
vector and a vector in the $x_1x_2$-plane. For a fixed value of $\theta$ this vector in the
$x_1x_2$-plane is fixed. As $\theta$ goes around the circle this vector rotates around once.

\begin{figure}[ht]
  \relabelbox \small {\epsfxsize=3in\centerline{\epsfbox{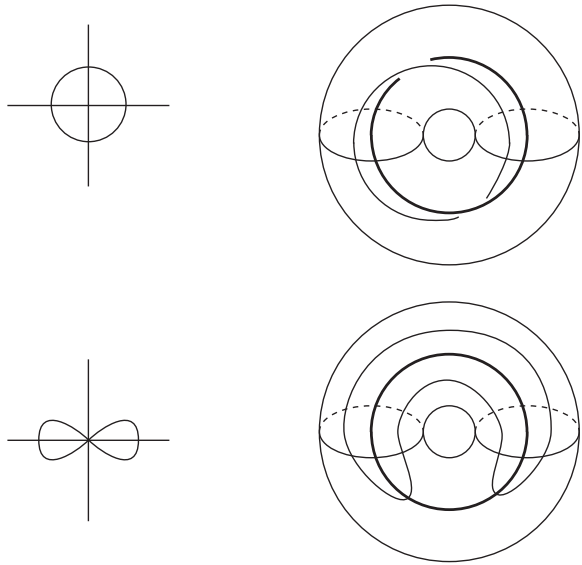}}}
  \endrelabelbox
        \caption{On the left hand side are two plane curves. On the right hand side is their
          corresponding Legendrian circles.}
        \label{fig:R2ex}
\end{figure}

Now given an oriented immersed curve $N=\gamma$ in $\R^2$ the Legendrian $L_N$ is simply the graph
of the ``twisted'' Gauss map. That is, if we let $R(\theta)= \theta+\frac{\pi}{2}$ and $g:N\to S^1$
be the Gauss map, then the twisted Gauss map is $R\circ g.$ (Note we are identifying $T\R^2$ and
$T^*\R^2$ using the flat metric on $\R^2.$) See Figure~\ref{fig:R2ex} for some examples.

Note that if $\gamma_0$ and $\gamma_1$ are two regular homotopic curves in $\R^2$ then any regular
homotopy $\gamma_t$ between $\gamma_1$ and $\gamma_2$ gives a map $C_{\gamma_t}\colon
S^1\times[0,1]\to W\times[0,1]$ such that $C(\bullet,j)$ parameterizes $L_{\gamma_j}$, $j=0,1$. It
is straightforward to check that the self intersection number of this cylinder is independent of the
particular regular homotopy chosen and that any dangerous self-tangency moment of a regular homotopy
contributes $\pm 1$ to the self intersection number. Thus if the algebraic self intersection number
of $C_{\gamma_t}$ is non-zero for some particular choice of regular homotopy then the same holds
true for all regular homotpies and we may conclude the necessity of self-tangencies. Twice the
self intersection number of $C_{\gamma_t}$ computes (one more than) the difference between a
``relative'' Thurston-Benniquin type invariant of $L(\gamma_0)$ and $L(\gamma_1).$ The
Thurston-Bennequin invariant of a Legendrian $S^1$ is a well-known invariant of Legendrian knots
in 3-dimensions. Thus we see shadows of contact geometry giving invariants of plane curves.

\subsection{General Co-dimension One immersions}\label{sec:codim1ex}
\begin{figure}[ht]
  \relabelbox \small {\epsfxsize=2in\centerline{\epsfbox{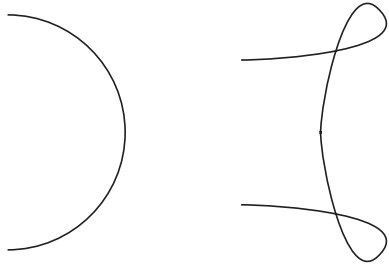}}}
  \endrelabelbox
        \caption{Two arcs in the $xy$-plane. The end points of each arc is on the $y$-axis.
          The spheres $S_1$ and $S_2$ are obtained form spinning the left and right arcs about
          the $y$-axis, respectively.}
        \label{fig:codim1}
\end{figure}

In \cite{Goryunov}, Goryunov studied immersions of surfaces in $\R^3$ (and more general codimension
one immersions) from a perspective similar to that described above. His methods allows to conclude
that there must be self tangencies in certain regular homotopies. One may also use contact homology
(defined below) in a similar way as we illustrate with the following example. Consider the two
immersions in Figure~\ref{fig:codim1}. Let $L_1$ and $L_2$ be the lifts of the fronts the left and
right ones respectively. While no ``classical invariants'' in contact geometry distinguish these two
Legendrian $S^2$'s contact homology will distinguish them. See Section~\ref{sec:compcd1ex} below for
this computation.

\subsection{Knots in $R^3$}
A knot $K$ is an embedded $S^1$ in $\R^3,$ so our ambient manifold is $M=\R^3.$ Thus the contact
manifold of interest is the unit cotangent bundle $W=\R^3\times S^2.$ The contact from on $W$ is
\[
\alpha= (q_1\, dp_1+q_2\, dp_2 + q_3\, dp_3)|_{W}.
\]
The knot $K$ is an embedded $S^1$ so the associated Legendrian is $L_K=S^1\times S^1.$ That is to a
knot $K$ in $\R^3$ we can associated a Legendrian torus $L_K=T^2$ in $W=\R^3\times S^2.$

What are the classical invariants of $L_K$ in $W$? The most obvious one is the homology class of
$L_K.$ This is determined by the degree of the map $T^2\to S^2$ obtained by the inclusion of $L_K$
into $W$ followed by projection of $W=\R^3\times S^2\to S^2.$ Note that this degree is simply the
normal degree of an infinitesimal torus around the knot and hence it equals $0$ by the Gauss-Bonnet
theorem. Legendrian submanifolds in any dimension have a Thurston-Bennequin invariant. If a
Legendrian submanifold in the contact manifold $W$ is zero-homologous then this is always equal to
half the Euler characteristic of the Legendrian submanifold. So for $L_K$ we always get 0 regardless
of the knot $K$ we started with. Finally there is the rotation class of a Legendrian submanifold
\[
r(L_K)\co H_1(L)\to \Z.
\]
This is defined as follows: given $\gamma$ a curve representing an element of $H_1(L)$ we can find a
surface $\Sigma$ in $\R^3\times S^2$ bounded by $\gamma.$ We can find a symplectic trivialization of
$\xi$ over $\Sigma,$ recall the two form $d\alpha$ gives $\xi$ a symplectic structure. Thus at each
point $x\in \gamma$ we can think of the Lagrangian plane $T_x L\subset \xi_x$ as a Lagrangian plane
in $\R^4$ with its standard symplectic structure. So $\gamma$ provides a loop of Lagrangian planes
in $\R^4.$ It is well known that $\mathcal{L}(2)$, the Grassmanian of Lagrangian planes in $\R^4,$
satisfies
\[
\pi_1(\mathcal{L}(2))=\Z.
\]
Thus to $\gamma$ we get an integer $r(L_K)(\gamma).$ In this situation one may show that this
integer is independent of $\Sigma$ and the trivialization of $\xi|_\Sigma.$ It can be shown that
\[
r(L_K)\equiv 0.
\]

Thus the classical invariants of a Legendrian submanifold give no interesting invariants of the knot
$K.$ This is one of the reasons the geometric construction we are considering here was not used to
study knots in $\R^3.$ We will see below that contact homology allows one to get interesting
invariants of knots in $\R^3$ from $L_K.$ Before we begin that part of the story we need to find
another way of thinking about the contact manifold $(W,\xi).$ We do this in a more general context
next.

\subsection{Submanifolds of Euclidean Space}
Here we consider any embedded $N$ in $\R^m.$ So our contact manifold is $W=\R^m\times S^{m-1}$ with
the contact form $\alpha = \sum_{i=1}^m p_i\, dq_i.$ Moreover our Legendrian $L_N$ is a
$(m-n-1)$-sphere bundle over $N.$ We wish to describe a new way of thinking about $(W,\alpha).$ To
this end consider the \dfn{1-jet space} of $S^{m-1}.$ The 1-jet space of $S^{m-1}$ is
\[
J^1(S^{m-1})= T^*S^{m-1} \times \R.
\]
Note $J^1(S^{m-1})$ is clearly a bundle over $S^{m-1}.$ Given any function $f:S^{n-1}\to \R$ its
1-jet is the section $j^1{f}(x)=(df(x),f(x))$ of $J^1(S^{m-1}).$ On $T^*S^{m-1}$ we have the
Louisville form $\lambda_{S^{m-1}}.$ The form $d\lambda_{S^{m-1}}$ is a symplectic structure on
$T^*S^{m-1}.$ We also have the contact from $\alpha=dz-\lambda_{S^{m-1}}$ on $J^1(S^{m-1}),$ where
$z$ is the coordinate on $\R.$

We will denote a point in $W=\R^m\times S^{m-1}$ by $({\bf q},{\bf p})$ where ${\bf q}=(q_1,\ldots,
q_m)$ and ${\bf p}=(p_1,\ldots, p_m)$ is a vector of unit length (again we identify the tangent and
cotangent bundles using the flat metric on $\R^m$). So ${\bf p}$ is a point on the unit sphere in
$\R^m,$ $\langle {\bf q}, {\bf p}\rangle$ is the part of ${\bf q}$ that is normal to the sphere at
${\bf p},$ and ${\bf q} - \langle {\bf q}, {\bf p}\rangle{\bf p}$ is the part of ${\bf q}$ that is
tangent to the sphere at ${\bf p}.$ We can define the map
\[
\Psi: W_{\R^m} \to J^1(S^{m-1})
\]
by
\[
\Psi({\bf p},{\bf q}) = ({\bf p}, {\bf q} - \langle {\bf q}, {\bf p}\rangle{\bf p}, \langle {\bf q}, {\bf p}\rangle).
\]
One may easily check that $\Psi$ is a diffeomorphism (in fact it is a bundle map covering the
identity on $S^{m-1}$ and is how one sees that $TS^{m-1}$ is stably trivial). Moreover $\Psi$ is a
contactomorphism, that is it takes the natural contact structure on $W$ to the natural contact
structure on $J^1(S^{m-1}).$

Thus when studying the Legendrian $L_N$ in $W$ we can instead study $L_N$ as a Legendrian
submanifold of $J^1(S^{m-1}).$ This has several advantages. Specifically there are two very helpful
projections of $J^1(S^{m-1}).$ First there is the front projection
\[
F: J^1(S^{m-1})\to S^{m-1}\times \R
\]
that just projects out the cotangent directions. If $(q_1,\ldots, q_{m-1})$ are local coordinates on
$S^{m-1}$ then $\lambda_{S^{m-1}}=\sum p_i\, dq_i$ and $\alpha$ on $J^1(S^{m-1})$ can be written
$dz-\sum p_i\, dq_i.$ In these local coordinates $F$ projects out the $p_i$-coordinates. Thus if we
consider the projection $F(L)$ of a Legendrian $L$ we see that the $p_i$-coordinates can be
recovered by looking at the slope of the tangent space to $F(L)$ (since $\alpha|_L=0$). Thus the
front projection allows us to reduce the dimension of the ambient space significantly (making it
easier to picture $L$) without loosing any information about $L$ in the total space! The second
projection is called the Lagrangian projection,
\[
\Pi: J^1(S^{m-1})\to T^*S^{m-1}
\]
projects out the $\R$ factor. It is easy to check that if $L$ is Legendrian in $J^1{S^{m-1}}$ then
$\Pi(L)$ is Lagrangian in $T^*S^{m-1}.$ It is also not hard to show that $L$ can be recovered (up to
translation in the $\R$ factor) from $\Pi(L),$ see \cite{ees1}. This Lagrangian projection will be
most useful to us in studying {\em contact homology}.

\section{Legendrian Contact Homology}
Here we define contact homology differential graded algebra (DGA) associated to a Legendrian
submanifold $L$ in a 1-jet space $J^1(S^{m-1}).$ (The definition easily generalizes to $J^1(M),$ or
indeed any exact symplectic manifold cross $\R,$ but to simplify the discussion we only discuss
$J^1(S^{m-1}).$) As mentioned above Legendrian contact homology is a small part of the newly defined
Symplectic Field Theory (SFT) of Eliashberg, Givental and Hofer \cite{egh}. While the analytic
underpinnings of the general theory of SFT are still being worked out, the foundations of Legendrian
contact homology in jet spaces has already been worked out in \cite{ees4}. Here we briefly describe
this theory.

We begin by fixing an almost complex structure $J$ on $T^*S^{m-1}.$ That is $J$ is a bundle
isomorphism of $T(T^*S^{m-1})\to T^*S^{m-1}$ such that $J^2=-id_{T^*S^{m-1}}.$  The contact homology
of $L$ will be a differential graded algebra (DGA).\hfill\break
{\bf The Algebra.} Denote the double points of the Lagrangian projection, $\Pi(L),$ by
$\mathcal{C}.$ We assume $\CC$ is a finite set of transverse double points. Let $\mathcal{A}$ be the
free associative unital algebra over $\Z_2$ generated by $\mathcal{C}.$\hfill\break
{\bf The Grading.} To each crossing $c\in\CC$ there are two points $c^+$ and $c^-$ in $L\subset
J^1(S^{m-1})$ that project to $c.$ We denote by $c^+$ the point with larger $z$-coordinate. Choose a
map $\gamma_c:[0,1]\to L$ that parametrizes an arc running from $c^+$ to $c^-.$ (Note there could be
more that one path.) For each point $\gamma(\theta)\in L$ we have a Lagrangian plane
$d\pi_{\gamma(\theta)}(T_{\gamma(\theta)}L)$ in $T^*S^{m-1}.$ Thus $\gamma$ give us a path
$\widehat{\gamma}$ in the bundle of Lagrangian subspaces $\text{Lag}(T(T^\ast S^{m-1})).$ Since $c$
is a transverse double point $\widehat{\gamma}(0)$ is transverse to $\widehat{\gamma}(1).$ Thus we
can find a complex structure $J'$ on $T_{\gamma(0)}T^*S^{m-1}$ (unrelated to $J$!) that (1) induces
the same orientation on $T_{\gamma(0) T^\ast S^{m-1}}$ as the almost complex structure $J$ and (2)
$J'(\widehat{\gamma}(1))=\widehat{\gamma}(0).$ Now set $\overline{\gamma}(\theta)=
e^{J'\theta}\widehat{\gamma}(1).$ Note that $\widehat{\gamma}$ followed by $\overline{\gamma}$ is a
closed loop in $\text{Lag}(T(T^\ast S^{m-1})).$ Moreover $\pi_1(\text{Lag}(T(T^\ast S^{m-1}))=\Z$
($m\ge 3$) thus to this closed loop we get an integer, $cz(c),$ the Conley-Zehnder or Maslov index
of $c.$ (This index can be computed by counting intersections between $\widehat\gamma$ and the
section $T^\ast S^{m-1}\to\text{Lag}(T(T^\ast S^{m-1}))$ which associates to each point the vertical
tangent space at that point.) The grading on $c$ is \[|c|=cz(c)-1.\] We note that $c$ in general
depends on the path $\gamma$ with which we started. To take care of this ambiguity note that to each
circle immersed in $L$ we get an integer through the above procedure. Let $n$ be the greatest common
divisor of all these integers. It is easy to convince oneself that $|c|$ is well defined modulo $n.$
\hfill\break
{\bf The Differential.} We will define the differential $\partial$ on $\A$ by defining it on the
generators of $\A$ and then extending by the signed Leibniz rule:
\[
\partial ab= (\partial a)b+(-1)^{|a|}a\partial b.
\]
Let $a\in \CC$ be a generator of $\A$ and let $b_1\ldots b_k$ be a word in the ``letters'' \CC. Let
$P_{k+1}$ be a $k+1$ sided polygon in $\C$ with vertices labeled counterclockwise $v_0,\ldots, v_k.$
We will be consider maps $u:(P_{k+1},\partial P_{k+1})\to (T^*S^{m-1}, \Pi(L))$ such that
$u|_{\partial P_{k+1}\setminus \{v_i\}}$ lifts to a map to $L\subset J^1(S^{m-1}).$ Call a vertex
$v_i$ mapping to the double point $c$ is \dfn{positive} (resp. \dfn{negative}) if the lift of the
arc just clockwise of $v_i$ in $\partial P_{k+1}$ lifts to an arc approaching $c^+$ (resp. $c^-$)
and the arc just counterclockwise of $v_i$ lifts to an arc approaching $c^-$ (resp. $c^+$), where
$c^\pm$ are as in the definition of grading. Set
\[
\M^a_{b_1\ldots b_k}=
\{u:(P_{k+1},\partial P_{k+1})\to (T^*S^{m-1}, \Pi(L)) \text{ satisfying 1.-- 4. below}\}/ \sim
\]
where $\sim$ is holomorphic reparameterization (which is relevant if $0\le k\le 1$) and the
conditions are
\begin{itemize}
\item[1.] $u|_{\partial P_{k+1}\setminus \{v_i\}}$ lifts to a map to $L\subset \R^{2n+1}$
\item[2.] $u(v_0)=a$ and $v_0$ is positive.
\item[3.] $u(v_i)=b_i, i=1,\ldots, k,$ and $v_i$ is negative.
\item[4.] $u$ is $J$-holomorphic.
\end{itemize}
We can now define
\[
\partial a=\sum_{b_1,\ldots, b_k} (\#_2\M) b_1b_2\ldots b_k,
\]
where the sum is taken over all words in the letters \CC for which $\dim(\M^a_{b_1\ldots b_n})=0$
and $\#_2$ denotes the modulo two count of elements in $\M.$
\begin{thm}[Ekholm, Etnyre and Sullivan, \cite{ees1, ees2, ees4}]
With the notation above:
\begin{enumerate}
\item The map $\partial$ is a well defined differential that reduces the grading by 1.
\item The stable tame isomorphism class of $(\A,\partial)$ is an invariant of $L.$
\item The homology of $(\A,\partial)$ is an invariant of $L.$
\end{enumerate}
\end{thm}
\section{Knot Invariants}
In this section we will indicate how to compute the Legendrian contact homology for the Legendrian
$L_K$ in $J^1(S^2)$ associated to a knot $K$ in $\R$ and thus have a (new?) invariants of $K.$ More
specifically we will explicitly compute the contact homology associated to an unknot and then show
how to compute the contact homology of a general knot in terms of this and a braid representation of
the knot. Lenny Ng has carried out this analysis in \cite{Ng1} and predicted what the contact
homology of a knot should be in terms of a braid representing the knot. He then proceeded to show
that this gives an invariant of a knot by showing the predicted contact homology DGA is independent
of the specific braid used to represent the knot. After the work in \cite{ees4} it is known that the
contact homology of $L_K$ is a well-defined invariant of $K.$ It is an ongoing research project to
verify that Ng's DGA is indeed the contact homology DGA.

\subsection{The Unknot}\label{sec:uk}
\begin{figure}[ht]
  \relabelbox \small {\epsfxsize=3in\centerline{\epsfbox{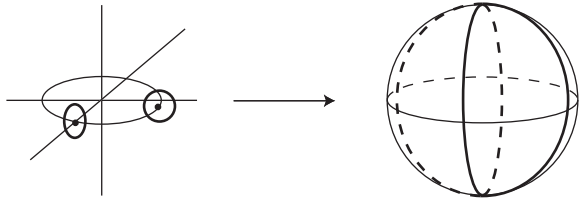}}} \endrelabelbox
        \caption{Unit circle in $xy$-plane with two unit normal circles indicated (left). The
          image of these unit circles in the unit cotangent sphere (right).}
        \label{fig:IS2}
\end{figure}

\begin{figure}[ht]
  \relabelbox \small {\epsfxsize=3in\centerline{\epsfbox{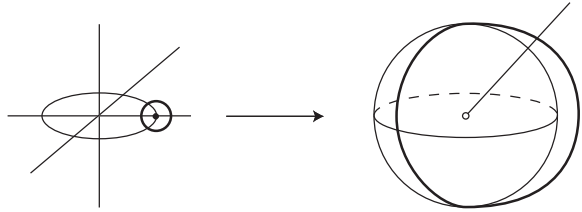}}}
  \endrelabelbox
        \caption{Unit circle in $xy$-plane with one unit normal circle indicated (left). The image of these unit circle
        in $S^2\times\R=\R^3\setminus \{(0,0,0)\}$ (right).}
        \label{fig:IF}
\end{figure}

\begin{figure}[ht]
  \relabelbox \small {\epsfxsize=1.2in\centerline{\epsfbox{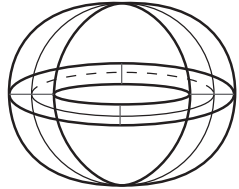}}} \endrelabelbox
        \caption{The front projection of the unknot. The four arcs connect points on the front that
          project to double points in the Lagrangian projection (after the original $S^1$ is perturbed
          to an ellipse).}
        \label{fig:UF}
\end{figure}

We begin with the simplest knot the round unknot $U$ sitting in the $xy$-plane in $\R^3.$ To
understand $L_U$ we draw its front projection and we start drawing the front projection by
considering the image of $L_U$ in $S^2.$ To understand this consider Figure~\ref{fig:IS2}.
Recall $S^2$ is thought of as the unit 2-sphere in the cotangent space (which we identify with the
tangent space) at each point of $\R^3.$ Thus to each point on $U$ the unit (co)normal bundle is just
the unit $S^1$ in the plane orthogonal to $U$ at that point. Figure~\ref{fig:IS2} draws a few of
these $S^1$'s. At each point of $U$ we see the $S^1$ we get is a great circle on $S^2$ and goes
trough the north and south poles. As we traverse $U$ these great spheres rotate around $S^2.$ We now
draw the front projection. This will be a $T^2$ in $S^2\times \R.$ We will think of $S^2\times \R$
as $\R^3$ with the origin removed. That is $S^2\times\{0\}$ is the unit $S^2$ in $\R^3$ and
$S^2\times\{p\}$ for $p<0$ are concentric spheres inside the unit sphere and for $p>0$ they are
concentric spheres outside the unit sphere. In Figure~\ref{fig:IF} we draw the image of one of the
great circle $S^1$'s in $S^2\times \R.$
Recall the \R factor is simply the normal component of the position vector (in $U$) at the point in
the unit normal $S^1.$ Clearly above the north and south poles this normal component is 0 and along
the equator the magnitude of the normal component is maximal and of opposite sign for the two points
intersecting the equator. Thus one gets the circle pictured in Figure~\ref{fig:IF}. Since the
picture is symmetric we get the entire front projection by simply rotating this $S^1$ about the axis
through the north and south pole. See Figure~\ref{fig:UF}.

\begin{figure}[ht]
  \relabelbox \small {\epsfxsize=1.2in\centerline{\epsfbox{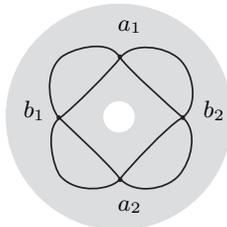}}}
  \relabel{1}{$b_1$}
  \relabel{2}{$b_2$}
  \relabel{3}{$a_2$}
  \relabel{4}{$a_1$}
  \endrelabelbox
        \caption{The Lagrangian projection of the Legendrian above the equator in $S^2.$}
        \label{fig:2ds}
\end{figure}

We now want to compute the contact homology of $L_U.$ So where are the double points in the
Lagrangian projection? In the front projection these correspond to pairs of points in $L_U$ with the
same $S^2$ coordinate and parallel tangent planes. Thus it is easy to see we have an $S^1$ of double
points over the equator of $S^2.$ This is a very degenerate situation! To fix this we simply isotop
$U$ into an ellipse in the $xy$-plane. It is clear that $L_U$ looks very similar, but now there are
only four double points. See Figure~\ref{fig:UF}. We denote the double points $a_1,a_2, b_1$ and
$b_2$ as indicated in the figure. One can check that $|a_i|=1$ and $|b_i|=2.$ Thus we have $\partial
a_i=0.$ Actually it is possible that $\partial a_i$ is 1 if there is an odd number of (rigid)
holomorphic disk with one positive puncture at $a_i.$ It is not obvious that there are any such
disks, but in fact there are four for each $a_i$! Two running up towards the north pole and two
running down towards the south pole. One may explicitly find these disks in $T^*S^2$ and show these
are the only possible disks. This computation also follows from the ideas in Section~\ref{cft}
below. In any event $\partial a_i=0.$ To understand the differential of the $b_i$ note that we can
think of $S^1$ as the equator of $S^2$ and thus $J^1S^1\subset J^1S^2.$ Now consider the Lagrangian
projection $\Pi(L_U)$ and intersect it with $T^*S^1\subset T^*S^2.$ We can arrange that $T^*S^1$ is
a complex subset of $T^*S^2$ and thus if we see holomorphic disks in $T^*S^1$ they will also be
holomorphic disks in $T^*S^2.$ We show $T^*S^1$ in Figure~\ref{fig:2ds}.
Clearly there are four holomorphic disks with boundary on $L_U\cap T^*S^1.$ It is not hard to show
these contribute to the differential of the $b_i$'s and that there are no other holomorphic disks
that do. So we have
\[
\partial a_i=0, \quad \quad \partial b_i=a_1+ a_2.
\]
Thus we completely understand the DGA of the Legendrian associated to an unknot! Unfortunately it is
difficult to directly compute the DGA for any other knot, but due to a clever ``localization'' idea
of Eliashberg we can still compute the DGA using braid theory.

\subsection{Braids and gradient flow trees}\label{sec:braid}
We now consider an arbitrary knot $K$ in $\R^3$ it is a well known result of
Alexander~\cite{Alexander23} that $K$ may be braided about the unknot $U$. That is we can fix a
tubular neighborhood $N=U\times D^2=S^1\times D^2$ and then isotop $K$ into $N$ so that $K$ is
transverse to $\{x\}\times D^2$ for all $x$ in $S^1.$ Once we have done this, we can shrink the
radius of the disks $D^2$ in $N.$ As we do this $K$ is approaching $U$ (as a $n$ fold cover if $K$
wraps $n$ times around $U$) and the tangent vectors of $K$ are getting very close to the tangent
vectors of $U.$ This says that as the radii of the disk $D^2$ shrink the Legendrian torus $L_K$ is
getting arbitrarily close to $L_U.$

A theorem of Weinstein \cite{Weinstein71} says that any Lagrangian $L$ in a symplectic manifold
$(X,\omega)$ has an neighborhood $N_L$ that is symplectomorphic to a neighborhood $N_0$ of the
0-section in $T^*L$ (with the standard symplectic form on it). Thus, returning to our braided knot
$K,$ we can think of $L_K$ as sitting inside $N_0$ and, in particular, inside $T^*T^2.$ It is not
too hard to see that $L_K$ is an exact Lagrangian in $T^*T^2$ and thus lifts to a Legendrian in
$J^1(T^2).$ So to a braided knot $K$ we get a Legendrian torus in $J^1(T^2).$ Isotopies of $K$ as a
braid (not a knot!) yield Legendrian isotopies of the Legendrian $L_K$ in $J^1(T^2).$ Thus studying
the Legendrian contact homology of $L_K$ in $J^1(T^2)$ will give us an invariant of braids. When
thinking of $L_K$ as a Legendrian in $J^1(T^2)$ we will denote it as $L_B,$ the $B$ means braid. In
fact we should think of $B$ as a particular representation of $K$ as a braid.

It is actually much easier to study the contact homology of $L_B$ than that of $L_K.$ We first
give a description of $L_B$ as the 1-jet of a ``multi-function''. To do this we will think of $B$ as
a multi function in $N=S^1\times D^2.$ That is we can find $n$ functions $f_1,\ldots, f_n\co
[0,1]\to D^2$ such that $f_i(t)\not= f_j(t)$ for $i\not=j$ and such that the sets $\{f_i(0)\}$ and
$\{f_i(1)\}$ are equal. Moreover gluing $\{0\}\times D^2$ to $\{1\}\times D^2$ by the identity map
gives $N=S^1\times D^2$ and the image of the graphs of the $f_i$ in $N$ equal $B.$ Said another way
we can think of $F(t)=(f_1(t),\ldots, f_n(t))$ as a function from $S^1$ to the set of distinct
$n$-tuples of points in $D^2$ and the ``graph'' of $F$ is $B.$ (Note if $n=1$ then we are talking
about an honest function.) We can describe $L_B$ as the 1-jet of the ``multi-function'' defined by
\[
g_i\co [0,1]\times[0,1]\to \R\co (t,\theta)\mapsto \langle f_i(t), (\cos 2\pi \theta, \sin 2\pi
\theta)\rangle.
\]
So the front projection $F(L_B)$ is given by the graphs of these functions in $T^2\times \R.$

We first look for the double point of the Lagrangian projection. From our discussion of the unknot
above we know the double points of $\Pi(L_B)$ will be at pairs of points in $F(L_B)$ with the same
$T^2$ coordinates and with parallel tangent planes. To say the tangent planes are parallel is
equivalent to saying $dg_i=dg_j,$ for some $i,j.$ Or if we set
\[
g_{ij}(t,\theta)=g_i(t,\theta)-g_j(t,\theta)=\langle f_i(t)-f_j(t), (\cos 2\pi \theta, \sin 2\pi
\theta)\rangle
\]
then the tangent planes are parallel where $dg_{ij}=0.$ So critical points of the difference
function $g_{ij}$ give double points of $\Pi(L_B).$ It is easy to see that these critical points
occur where $|f_i-f_j|$ has a critical point. One may arrange that $|f_i-f_j|:S^1\to \R$ has exactly
one minimum and one maximum which we denote $a_{ij}$ and $b_{ij},$ respectively. The grading on
these double points is related to their Morse index and can be computed to be $|a_{ij}|=0$ and
$|b_{ij}|=1.$

Moving to the differential we clearly have $\partial a_{ij}=0$ since the differential reduces the
grading by one. The computation of the differential on the $b_{ij}$'s is much more complicated. We
begin with a ``simple'' toy situation first studied by Floer. Suppose $L_B$ has two components each
of which is simply the graph of the 1-jet of a function. Thus $L_B$ is described by two functions
$g_1:T^2\to \R$ and $g_2:T^2\to \R.$ (Note because there are two components one cannot grade the DGA
as we did above. We will ignore this subtlety here as we are just considering this situation to get
an idea what the double points of the Lagrangian projection are and what the holomorphic disks are.)
Floer showed that if $g_{12}=g_1-g_2$ is sufficiently small there is a one to one correspondence
between the gradient flow lines of $g_{12}$ and holomorphic disk in $T^*T^2$ connecting double
points. Thus to compute the differential of a double point/critical point $a$ one just counts
``rigid gradient flow lines''. We call a flow line rigid if there are no flow lines with the same
endpoints near the given one. So
\[
\partial a = \sum_{b} \#_2 \mathcal{F}^a_b b
\]
where $\mathcal{F}^a_b$ are the flow lines from $a$ to $b$ and the sum is taken over $b$ for which
$\mathcal{F}^a_b$ consists of rigid flow lines.

\begin{figure}[ht]
  \relabelbox \small {\epsfxsize=3.5in\centerline{\epsfbox{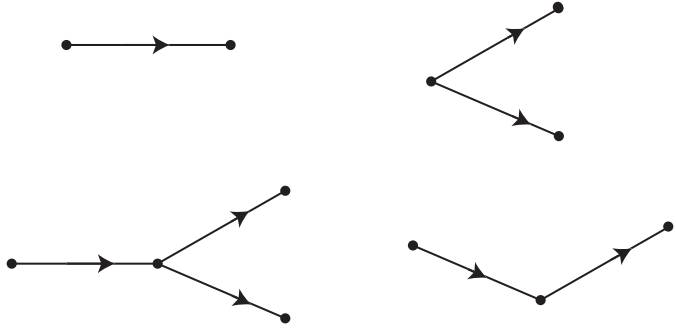}}}
  \endrelabelbox
        \caption{Various types of flow trees as seen on $T^2.$}
        \label{fig:ft}
\end{figure}

\begin{figure}[ht]
  \relabelbox \small {\epsfxsize=4in\centerline{\epsfbox{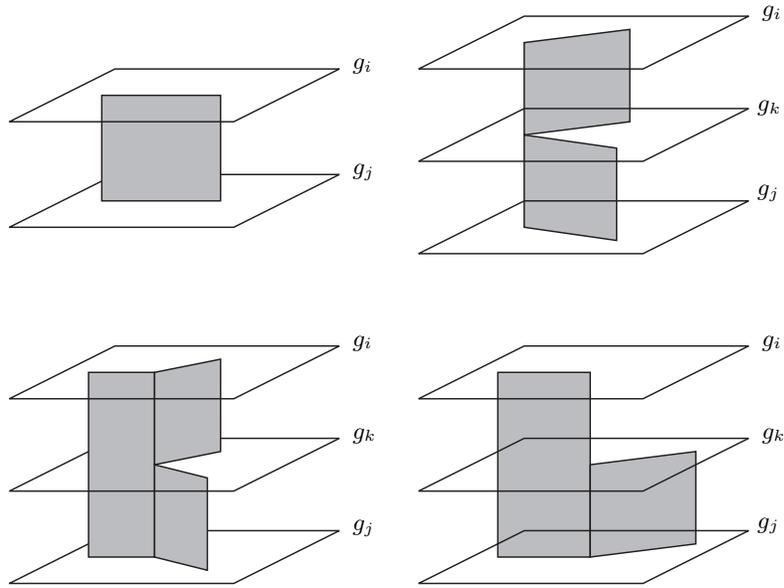}}}
  \relabel{1}{$g_i$}
  \relabel{2}{$g_j$}
  \relabel{3}{$g_i$}
  \relabel{4}{$g_k$}
  \relabel{5}{$g_j$}
  \relabel{6}{$g_i$}
  \relabel{7}{$g_k$}
  \relabel{8}{$g_j$}
  \relabel{9}{$g_i$}
  \relabel{a}{$g_k$}
  \relabel{b}{$g_j$}
  \endrelabelbox
        \caption{The disks lying above the flow trees in Figure~\ref{fig:ft}. On the bottom right we
          show only one type of flow tree associated to a ``non-source'' valence 2 vertex.}
        \label{fig:ftd}
\end{figure}

\begin{figure}[ht]
  \relabelbox \small {\epsfxsize=1.6in\centerline{\epsfbox{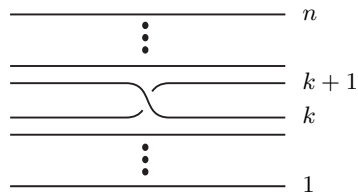}}}
  \relabel{1}{$1$}
  \relabel{k}{$k$}
  \relabel{p}{$k+1$}
  \relabel{n}{$n$}
  \endrelabelbox
        \caption{A braid with $n$ strands and a single crossing between the $k$ and $k+1$st strands.}
        \label{fig:bk}
\end{figure}

Having understood this simpler situation we now return to the contact homology of $L_B$ and in
particular the computation of the differential of the $b_{ij}$'s. When we just had two functions
involved we reduced counting holomorphic curves to counting gradient flow lines. In our current
situation where there will be more than two functions and we must count {\em gradient flow trees}.
In our situation a gradient flow tree is a tree $G$ and an embedding $f:G\to T^2$ where
\begin{enumerate}
\item $G$ has only 1, 2 and 3 valent verticies and oriented edges such that exactly one valence 1 or
 2 vertex having all the edge orientation(s) pointing away from it
 \item $f$ maps the 1 and 2 valent verticies of $G$ to critical points of the $g_{ij}$'s,
\item for each edge there is an $i,j$ such that $f$ map the edge to part of a flow line of $g_{ij}$
  so that the orientation points in the direction of the flow, and
\item each 3-valent vertex maps to a point that is not a critical point of any $g_{ij},$ one edge is
  pointing to the vertex and is labeled $g_{ij},$ the other two edges point away from the vertex and
  there is a $k$ such that $g_i>g_k>g_j$ near the vertex and the two two edges pointing away from
  the vertex are labeled $g_{ik}$ and $g_{kj}$
\item at a valance 2 vertex with both edges pointing away from it, the edges are labeled $g_{ik}$
and $g_{kj}$ where $g_i>g_k>g_j$ near the vertex and the vertex is mapped to a critical point of
$g_{ij}$
\item at any other valence 2 vertex one edge points towards it and is labeled $g_{ij},$ one points
away from it and is labeled $g_{ik}$ or $g_{kj}$ and the vertex is mapped to a critical point of
$g_{ki}$ or $g_{kj},$ respectively, where $g_i>g_k>g_j$ near the vertex.
\end{enumerate}
See Figure~\ref{fig:ft} for examples of ``flow trees''
and Figure~\ref{fig:ftd} for the ``disks lying above them'' in $T^2\times \R.$
We call the valence 1 and 2 verticies the corners of the tree. A flow tree is called rigid if there
are no near by flow trees with the same corners. Let $\mathcal{T}^a_{b_1\ldots b_l}$ be the set of
flow trees with a corner mapping to $a$ and reading counterclockwise around the boundary of the
``disk lying above the tree'' starting at $a$ we read the word $ab_1\ldots b_n.$ Following the ideas
in \cite{FO} one can show that the holomorphic curves in $T^*T^2$ correspond to flow trees (assuming
all the $g_{ij}$'s are sufficiently small). More precisely given the Legendrian $L$ expressed as the
graph of a multi function $\{g_i\}$ then let $L_\lambda$ be the Legendrian given by the multi
function $\{\lambda g_i\}.$ Clearly all the $L_\lambda$ are Legendrian isotopic.
\begin{thm}
For $\lambda$ sufficiently small there is a 1-1 correspondence between rigid holomorphic disks with
boundary on $\Pi(L_\lambda)$ and rigid gradient flow trees of the $g_{ij}=g_i-g_j.$
\end{thm}
Thus to compute the boundary of the double points $b_{ij}$ we can count the rigid gradient flow
trees for the $g_{ij}$'s. That is
\begin{equation}\label{gftboundary}
\partial b_{ij} = \sum \#_2 (\mathcal{F}^{b_{ij}}_{{a_{{i_1}{j_1}}}\ldots {a_{{i_l}{j_l}}}})
{a_{{i_1}{j_1}}}\ldots {a_{{i_l}{j_l}}},
\end{equation}
where the sum is over ${a_{{i_1}{j_1}}}\ldots {a_{{i_l}{j_l}}}$ such that the flow trees are rigid.
It is not too hard to compute this for a braid $\sigma_k$ shown in Figure~\ref{fig:bk}.
For a general braid Ng \cite{Ng1} computed
\begin{equation*}
    \partial b_{ij} = a_{ij} -  \phi_B(a_{ij}).
\end{equation*}
Where $\phi_B$ is defined as follows: for a generator $\sigma_k$ of the braid group we define
\begin{equation*}
\phi_{\sigma_k}:   \begin{cases}
        a_{ki} \mapsto -a_{k+1,i}-a_{k+1,k}a_{ki} & i \neq k, k+1 \\
        a_{ik} \mapsto -a_{i,k+1}-a_{ik}a_{k,k+1} & i\neq k, k+1\\
        a_{k+1,i} \mapsto a_{ki} & i\neq k, k+1\\
        a_{i,k+1} \mapsto a_{ik} & i\neq k, k+q\\
        a_{k,k+1} \mapsto a_{k+1,k} &\\
        a_{k+1,k} \mapsto a_{k,k+1} &\\
        a_{ij} \mapsto a_{ij} & i,j \neq k, k+1.
    \end{cases}
  \end{equation*}
If $B=\sigma_{i_1}\ldots \sigma_{i_l}$ then define $\phi_B=\phi_{\sigma_{i_l}}\circ \ldots\circ
\phi_{\sigma_{i_1}}.$ Thus we have computed the contact homology of the Legendrian $L_B$ associated
to the braid $B.$ This is an invariant of $B$! But to get an invariant of the knot $K$ that $B$
represents more work needs to be done. In particular we must glue the local answer obtained here
into $J^1(S^2).$ This will be discussed in the following sections.

\subsection{Interlude: Computation of the contact homology of immersed spheres in $\R^3$}
\label{sec:compcd1ex}
Recall from Section~\ref{sec:codim1ex} we had two Legendrian spheres $L_1$ and $L_2$ that we were
trying to distinguish. Note that neither $L_1$ nor $L_2$ has any Reeb chords so the contact homology
of both of them are simply $\Z_2$. However, let $L_\infty$ be a copy of the zero-section shifted
very high up in the $\R$-direction of $J^1(S^2)=T^\ast S^2\times\R$. Then if $L_1$ and $L_2$ would
be Legendrian isotopic so would the two links $L_\infty\cup L_1$ and $L_\infty\cup L_2$. After a small
perturbation of $L_1$ by a Morse function with two critical points $L_\infty\cup L_1$ has two Reeb
chords $a$ and $b$. Moreover using the Floer-Morse correspondence we see that $\pa a_1=0$ and $\pa b_1=0$. Thus the contact homology of $L_\infty\cup L_1$ equals $\Z_2\la a_1,b_1\ra$ where the grading difference between the two generators is $1$.
Also, it is straightforward to show that after a small perturbation, the link $L_\infty\cup L_2$ has two Reeb chords $a_2$ and $b_2$ but the grading difference in this case is $3$. Moreover, the lengths of the two Reeb chords can be made almost the same which implies that there no rigid holomorphic disks.
Therefore the contact homology of $L_\infty\cup L_2$ equals $\Z_2\la a_2,b_2\ra$ where the grading difference between the two generators is $3$. Since the
contact homology of the links differ we conclude that $L_1$ and $L_2$ are not Legendrian isotopic
and hence every regular homotopy between $L_1$ and $L_2$ must have a dangerous self tangency.

\subsection{Cusp Flow Trees}\label{cft}
In this subsection we return to calculating the boundary map for the DGA associated to the unknot.
Recall we have determined that the Legendrian torus $L_U$ associated to the elliptical unknot in the
$xy$-plane has four double points $a_1,a_2,b_1$ and $b_2$ in its Lagrangian projection. The gradings
are $|a_i|=1$ and $|b_i|=2.$ Moreover, we have seen that $\partial b_i=a_1+a_2.$ In
Section~\ref{sec:uk} we claimed that $\partial a_i=0$ but did not give a justification there. We now
indicate how this computation was made.

\begin{figure}[ht]
  \relabelbox \small {\epsfxsize=4in\centerline{\epsfbox{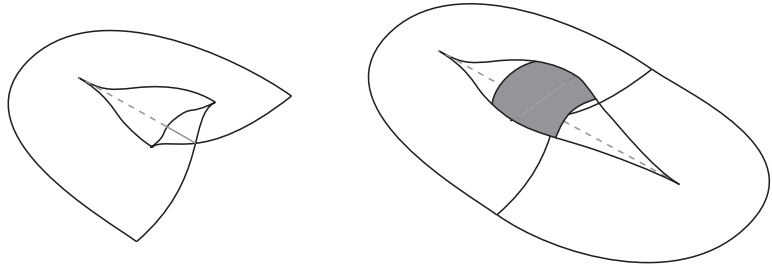}}}
  \endrelabelbox
        \caption{A swallow tail singularity, left. On the right is two swallow tail singularities
          glued together.}
        \label{fig:sw}
\end{figure}

\begin{figure}[ht]
  \relabelbox \small {\epsfxsize=4in\centerline{\epsfbox{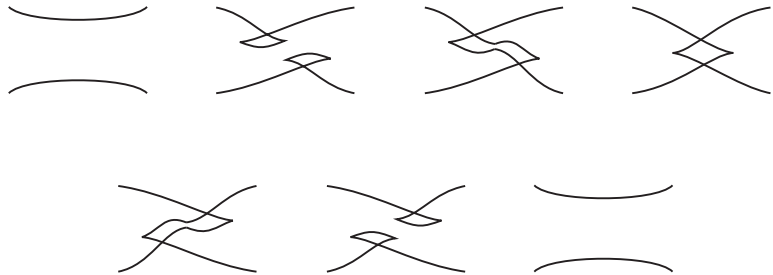}}}
  \endrelabelbox
        \caption{Slices of the Legendrian near the ``north'' and ``south pole''.}
        \label{fig:slice}
\end{figure}

To begin with we observe that the front projection of $L_U,$ shown in Figure~\ref{fig:UF}, is still
very degenerate above the north and south poles of $S^2.$ One can slightly perturb $L_U$ to make the
front generic or one can start with an unknot $U$ that is not in the $xy$-plane, but rather on the
graph of $z=\epsilon(x^2-y^2)$ for small $\epsilon.$ In the latter case the Legendrian $L_U$ will
look essentially the same as it did before except near the north and south pole. There quite a
change takes place. It is an interesting exercise to see the front changes as follows.
Figure~\ref{fig:sw} shows a ``swallow tail'' singularity on the left, it is ``doubled'' on the
right. Take two copies of the picture on the right and remover the grey shaded region from both
copies. Now turn one copy upside down and rotate it by $\frac\pi 2$ and glue the copies together
along the boundary of the grey shaded region. This produces a front of an annulus. See
Figure~\ref{fig:slice} for a different view of this annulus. In this figure slices of the front are
shown. If a small neighborhood of the singular point at the north pole is removed and replaces by
the annulus just constructed then the front will be generic near the north pole. Similarly we can do
the same at the south pole. The resulting front is a generic front of the Legendrian associated to
the unknot.

\begin{figure}[ht]
  \relabelbox \small {\epsfxsize=4in\centerline{\epsfbox{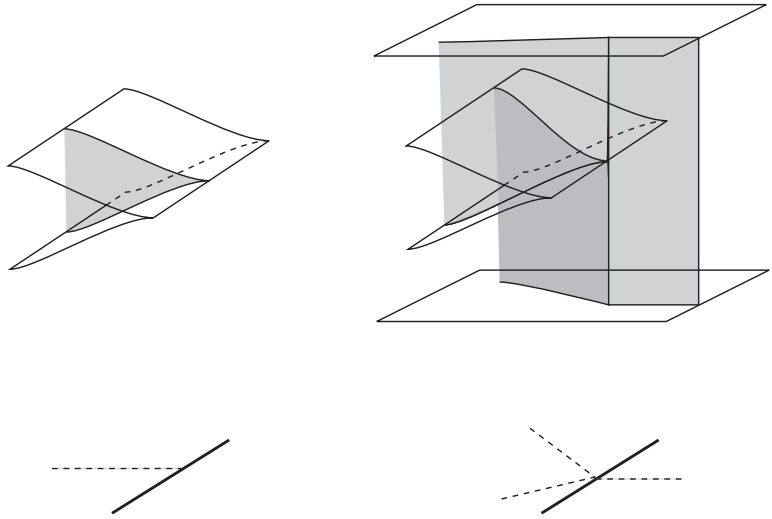}}}
  \endrelabelbox
        \caption{At the bottom of the figure are two ``cusped flow trees''. The thicker line is the
          projection of a cusp to the base manifold and the dotted lines are parts of a flow tree.
        on the left we see an edge ending in a cusp and on the right we see an edge split into two
      edges at a cusp. At the top of the figure are the ``disks lying above these flow trees''.}
        \label{fig:cft}
\end{figure}

This front can still be thought of the graph of a multi function, but now the number of values the
multi function takes changes over different regions of $S^2$ and some portions of the graph come
together along cusps. We will call such a multi function a ``cusped multi function''. We can still
consider gradient flow trees for this cusped multi function, but now two new phenomena can occur: a
tree can have an edge ending at a cusp or a tree can have a valance three vertex at a cusp where the
flow lines splits. We do not describe these ``cusped flow trees'' as thoroughly as we did the
gradient flow trees, but all main ideas are encoded in Figure~\ref{fig:cft}.
The boundary map for the DGA can now be defined as we did in Equation~\eqref{gftboundary} except now
we must count all rigid cusped flow trees and not just rigid gradient flow trees.

\begin{figure}[ht]
  \relabelbox \small {\epsfxsize=4in\centerline{\epsfbox{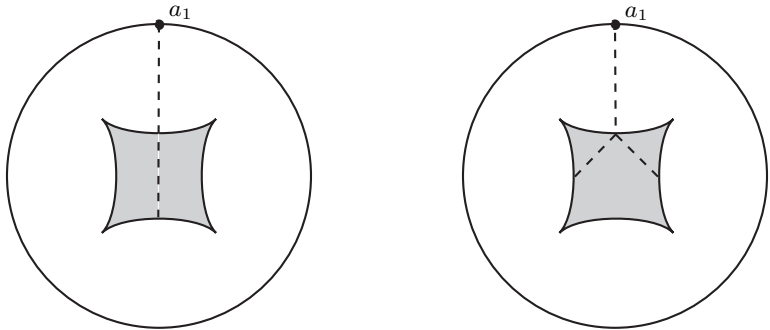}}}
  \relabel{a1}{$a_1$}
  \relabel{a2}{$a_1$}
  \endrelabelbox
        \caption{The northern hemisphere of $S^2.$ The multi function defining the graph has four values
          above the shaded region and two values elsewhere.}
        \label{fig:duk}
\end{figure}

With the above understood it is now a simple matter to compute $\partial a_i$ for the Legendrian
$L_U$ associated to the unknot. In particular, in Figure~\ref{fig:duk}
we see the northern hemisphere of $S^2$ on which the projection of the cusps of $L_U$ are drawn. The
two dotted trees indicate the only two rigid cusp flow trees emanating from $a_1.$ There is a
similar picture for the southern hemisphere. Thus as claimed in Section~\ref{sec:uk} there are four
holomorphic disks associated to $a_1$ so when counted modulo two we see $\partial a_i=0.$

\subsection{General Knots and Enhancements}\label{genpic}
Recall, when studying a general knot $K$ in $\R^3$ we braided it about the unknot $U.$ Thus the
Legendrian $L_K$ in $J^1(S^2)$ lay in a small neighborhood of $L_U$ which is contactomorphic to a
small neighborhood of the zero section in $J^1(T^2).$ In Section~\ref{sec:braid} we computed the
contact homology of $L_K$ in $J^1(T^2).$ We must now ``glue'' this back into $J^1(S^2).$ All the
double points of $L_K$ in $J^1(T^2)$ are still double points of $L_K\subset J^1(T^2)\subset
J^1(S^2),$ so $L_K$ in $J^1(S^2)$ has double points $a_{ij}, b_{ij}, 1\leq i,j,\leq n, i\not=j,$
where $n$ is the number of strands in the braid representing $K.$ Also recall $L_U$ has four double
points. Since $L_K$ is formed from $n$ copies of $L_U$ we see $L_K$ also has $4n^2$ more double
point $c_{ij}, d_{ij}, e_{ij}, f_{ij},$ where $c_{ij}$ and $d_{ij}$ correspond to the double points
for $L_U$ with grading 1, so they all have grading 1. Similarly the $e_{ij}$ and $f_{ij}$ all have
grading 2.

We have now identified all the generators of the contact homology of $L_K$ in $J^1(S^2).$ To
identifying the boundary map one must understand the holomorphic disks with boundary on $L_K$. Since
this requires a rather technical analysis of the situation we will merely give a sketch. The main
idea is to let $L_K$ approach $L_U$. In terms of the geometry of the braid this means that we push
the braid towards the circle along which it was braided. A rough statement of the result needed for
the computation is that for all $L_K$ sufficiently close to $L_U$ there is a $1-1$ correspondence
between rigid holomorphic disks with boundary on $L_K$ and the following objects.
\begin{itemize}
\item Gradient flow trees as described in Section \ref{sec:braid}.
\item Disks with boundary on $L_U$ with gradient trees emanating from their boundaries.
\end{itemize}
The idea is when the Legendrian $L_K$ ``degenerates'' on to $L_U$ all the holomorphic disk
associated to $L_K$ will break into disks in the ``localized'' situation considered in
Section~\ref{sec:braid} or into disks associated to $L_U$ or combinations of the two. In particular,
to compute the differential it is sufficient to know the gradient trees of $L_K\subset J^1(T^2)$ and
the moduli spaces of holomorphic disks with boundary on $L_U$.

This analysis is on going work of many people (Ekholm, Etnyre, Ng and Sullivan) and confirms
that Ng's predicted answer in \cite{Ng1} correspond to the actual holomorphic disk count. (Ng proved
that the stable tame isomorphism class of the DGA gives an invariant of knots in a combinatorial
fashion without reference to holomorphic disks.) Moreover, bringing orientations into the picture
(see \cite{ees3}) Ng has constructed an invariant of knots that should compute the contact homology
over the group ring $\Z[H_1(T^2)].$ He has shown that this is a very effective invariant. In
particular, one can extract the Alexander polynomial for the DGA and essentially the $A$-polynomial.
Thus the DGA will distinguish the unknot from all other knots! Moreover, the DGA is stronger than
these classical invariants as it can distinguish knots with the same Alexander polynomial and
$A$-polynomial. At the time of this writing it seems in the realm of possibility that this enhanced
DGA is a complete knot invariant!

The results needed to show Ng's DGA actually computes the contact homology DGA of $L_K$ hold in
greater generality. For example the correspondence between rigid cusped gradient flow trees and
rigid holomorphic disks hold in general for $1-$ and $2$-dimensional Legendrian submanifolds in
$1$-jet spaces and in the higher dimensions in the absence of front singularities other than cusp
edges.

This connection to holomorphic disks is very interesting as it gives a geometric meaning to Ng's
algebra. This could help establish properties of the DGA and, in addition, understanding Ng's DGA in
terms of contact homology would allow generalizations to other dimensions.


\end{document}